\documentclass[12pt]{article}

\usepackage{comment}
\usepackage{pgf,tikz}
\usepackage{pgfplots}
\usetikzlibrary{arrows.meta}

\usepackage{amsmath,amssymb,amsthm,cite}
\usepackage[colorlinks,citecolor=red,urlcolor=blue,linkcolor=blue]{hyperref}
\usepackage{mathrsfs,dsfont}
\usepackage{graphicx}
\usepackage{enumitem}

\newtheorem{nummer}{ }
\newtheorem{thm}[nummer]{\bf Theorem}
\newtheorem{prop}[nummer]{\bf Proposition}
\newtheorem{lem}[nummer]{\bf Lemma}
\newtheorem{cor}[nummer]{\bf Corollary}

\newtheorem{defi}[nummer]{\bf Definition}
\newtheorem{ques}{\bf Question}

\newtheorem{rmk}{\bf Remark}
\newcommand{\rmksty}{\rm}
\newtheorem*{alg}{\bf Algorithm}
\newcommand{\algsty}{\rm}

\newcommand{\ie} {\sl i.e.}
\newcommand{\eg} {\sl e.g.}

\newcommand{\pp}{pytha\-gorean pair}
\newcommand{\pps}{pytha\-gorean pairs}
\newcommand{\hl}{hl-pair}
\newcommand{\hls}{hl-pairs}
\newcommand{\hlppp}{pytha\-potent {\hl} of degree~$h$}
\newcommand{\hlpppd}[1]{pytha\-potent {\hl} of degree~$#1$}

\newcommand{\tor}[2]{\Z/#1\Z\times\Z/#2\Z}

\newcommand{\Q}{\mathds{Q}}

\newcommand{\Z}{\mathds{Z}}

\newcommand{\C}{\Gamma}
\newcommand{\Ccd}{\C_{c,d}}
\newcommand{\Ccdh}{\C_{c^h,d^h}}
\newcommand{\Ccah}{\C_{c^h,a^h}}
\newcommand{\Ccbh}{\C_{c^h,b^h}}

\makeatletter
\def\opargproof[#1]{\par\noindent {\bf #1 }}
 \makeatother

\begin{document}
\begin{center}
{\LARGE\bf Pairing Powers of Pytha\-gorean Pairs~II}

\medskip

{\small Lorenz Halbeisen}\\[1.2ex]
{\scriptsize Department of Mathematics, ETH Zentrum,
R\"amistrasse\;101, 8092 Z\"urich, Switzerland\\ lorenz.halbeisen@math.ethz.ch}\\[1.8ex]
{\small Norbert Hungerb\"uhler}\\[1.2ex]
{\scriptsize Department of Mathematics, ETH Zentrum,
R\"amistrasse\;101, 8092 Z\"urich, Switzerland\\ norbert.hungerbuehler@math.ethz.ch}\\[1.8ex]
{\small Arman Shamsi Zargar}\\[1.2ex]
{\scriptsize Department of Mathematics and Applications, University of Mohaghegh Ardabili,
Ardabil, Iran\\ zargar@uma.ac.ir}
\end{center}

\hspace{5ex}{\small{\it key-words\/}: pytha\-gorean triple, hypotenuse-leg pair, pytha\-potent pair, elliptic curve}

\hspace{5ex}{\small{\it 2020 Mathematics Subject
Classification\/}: {\bf 11D72} 11G05

\begin{abstract}\noindent
{\small Let $(a,b,c)$ be a pytha\-gorean triple, {\ie}, $a^2+b^2=c^2$. Then $(c,a)$ and $(c,b)$
are {\it hypotenuse-leg pairs\/} ({\it {\hls}\/}), {\ie}, pairs $(c,d)$ of positive integers
such that $c^2-d^2$ is a non-zero square. An {\hl} $(c,d)$ is called a 
{\it pythapotent {\rm{hl}}-pair of degree $h$\/} if there is a
{\pp} $(k,l)$, which is not a multiple of the {\pp} $(d,\sqrt{c^2-d^2})$, such that
$(c^hk,d^hl)$ is again an {\hl}. To each {\hl} $(c,d)$ and each $h\geq 1$ we assign the elliptic curve
$\Ccdh\colon y^2=x^3+(c^{2h}-d^{2h})x^2-c^{2h}d^{2h}x$ and show that its torsion group is
$\tor 2 2$ and that $(c,d)$ is a {\hlppp} if and only if $\Ccdh$ has positive rank over~$\Q$.
In this case there are infinitely many such {\pps} $(k,l)$, not multiples of each other.
Moreover, every {\hl} is a {\hlpppd 1} and a {\hlpppd 3}, and we give sufficient conditions
for degree~$2$. The two cases $(c,a)$ and $(c,b)$, where $a=m^2-n^2$ and $b=2mn$, are
worked out separately.}
\end{abstract}

\section{Introduction}

A {\it \pp\/} is a pair $(k,l)$ of positive integers such that $k^2+l^2$ is a square;
as usual we write $k^2+l^2=\Box$. In~\cite{HHZ} we called a {\pp} $(a,b)$ a
pytha\-potent pair of degree~$h$ if there is a {\pp} $(k,l)$, not a multiple of $(a,b)$,
such that $(a^hk,b^hl)$ is a {\pp} (the cases $h=1,2$ were previously studied in~\cite{JNT} and also in~\cite{HHZ}).
In the present paper we replace the sum of squares in the product pair by a difference of squares.

Let $(a,b,c)$ be a pytha\-gorean triple, {\ie}, $a^2+b^2=c^2$. Up to a common factor we may write
$$\begin{aligned}
	a &=m^2-n^2, & b &=2mn, & c &=m^2+n^2
\end{aligned}$$
with integers $m>n>0$. Then
$$\begin{aligned}
	c^2-a^2&=b^2=\Box &\text{and} && c^2-b^2&=a^2=\Box.
\end{aligned}$$
This leads to the following notion.

\begin{defi}
A pair $(c,d)$ of positive integers is a {\em hypotenuse-leg pair\/} (or {\em\hl\/} for short) if
$c^2-d^2$ is a non-zero square. For an {\hl} $(c,d)$ we always write $e:=\sqrt{c^2-d^2}$, so that
$(d,e,c)$ is a pytha\-gorean triple and $(d,e)$ is a {\pp}.
An {\hl} $(c,d)$ is called a {\em\hlppp\/} if there exists a {\pp} $(k,l)$, which is not a multiple of $(d,e)$,
such that $(c^hk,d^hl)$ is an {\hl}, {\ie},
\begin{equation}\label{eq:def}
	\begin{aligned}
		c^2-d^2 &=\Box, & k^2+l^2&=\Box, & \text{and} && (c^hk)^2-(d^hl)^2 &=\Box\neq 0.
	\end{aligned}
\end{equation}
\end{defi}

Here and in the following, two pairs $(k,l)$ and $(k',l')$ are called {\em multiples of each other\/} if
$(k,l)=\lambda\,(k',l')$ for some rational number $\lambda>0$. In this case we write $(k,l)\sim(k',l')$.
Since the last two conditions in~\eqref{eq:def} are homogeneous in $(k,l)$, we may always replace a solution $(k,l)$
by a multiple, in particular by a pair of relatively prime positive integers.

Thus, for a pytha\-gorean triple $(a,b,c)$ there are two cases:
\begin{itemize}
\item[] {\bf Case 1:} the {\hl} $(c,a)$, {\ie}, $d=a=m^2-n^2$ and $e=b=2mn$; we look for
{\pps} $(k,l)$ with $(c^hk)^2-(a^hl)^2=\Box$.
\item[] {\bf Case 2:} the {\hl} $(c,b)$, {\ie}, $d=b=2mn$ and $e=a=m^2-n^2$; we look for
{\pps} $(k,l)$ with $(c^hk)^2-(b^hl)^2=\Box$.
\end{itemize}

The general theory 
is the same for both cases and is
formulated for an arbitrary {\hl} $(c,d)$. For a fixed triple, however, the two cases lead to different
elliptic curves which may behave differently, and the explicit formulas in terms of $m$ and~$n$
differ.

Notice that the excluded pairs $(k,l)$ which are multiples of $(d,e)$ are trivial solutions for $h=1$: If $k=\lambda d$ and $l=\lambda e$, then
$$\begin{aligned}
	(ck)^2-(dl)^2 &=(\lambda cd)^2-(\lambda de)^2=
	\lambda^2d^2(c^2-e^2)=\lambda^2 d^4.
\end{aligned}$$

To the {\hl} $(c,d)$ and a positive integer $h$ we assign the elliptic curve
$$\begin{aligned}
	\Ccdh: \  y^2 &= x^3+(c^{2h}-d^{2h})x^2-c^{2h}d^{2h}x = x(x+c^{2h})(x-d^{2h}).
\end{aligned}$$
In the two cases above these are the curves $\Ccah$ and $\Ccbh$, respectively.
We will show that $\Ccdh$ has torsion group $\tor 2 2$ for every $h\geq 1$ 
and that $(c,d)$ is a {\hlppp} if and only if $\Ccdh$ has positive rank over~$\Q$. 
In this case, the points of infinite order on $\Ccdh$ produce infinitely many {\pps} $(k,l)$, not multiples of
each other, such that $(c^hk,d^hl)$ is an {\hl}. Moreover, every {\hl} is a {\hlpppd 1} and a
{\hlpppd 3} and for degree~$2$ we give sufficient conditions.

\noindent {\bf Example.} Consider the triple $(a,b,c)=(3,4,5)$, {\ie}, $m=2$, $n=1$.
The {\hl} $(5,4)$ (Case~2) is a {\hlpppd 4}: for the {\pp} $(k,l)=(181252,109725)$, with
$k^2+l^2=211877^2$, we have
$$\begin{aligned}
	(5^4\cdot 181252)^2-(4^4\cdot 109725)^2 &=109744700^2.
\end{aligned}$$
However, since the curve $\C_{5^2,4^2}$ has rank~$0$, the {\hl} $(5,4)$ is not a {\hlpppd 2},
whereas the {\hl} $(5,3)$ of the same triple (Case~1) is a {\hlpppd 2}:
$$\begin{aligned}
	(5^2\cdot 51)^2-(3^2\cdot 140)^2 &=195^2, & 51^2+140^2&=149^2.
\end{aligned}$$
We will see below how such examples can be generated systematically.

\begin{rmk}\label{rmk:B}
{\rmksty Equivalently, one may ask for an {\hl} $(l',k')$ such that the ``crossed'' product
$(c^hk',d^hl')$ is a {\pp}, {\ie},
$$\begin{aligned}
	c^2-d^2 &=\Box, & l'^2-k'^2=\Box, && \text{and} && (c^hk')^2+(d^hl')^2 &=\Box.
\end{aligned}$$
Indeed, if $(k,l)$ satisfies~\eqref{eq:def}, then $(l',k'):=(c^hk,\,d^hl)$ satisfies
$l'^2-k'^2=(c^hk)^2-(d^hl)^2=\Box$ and $(c^hk')^2+(d^hl')^2=c^{2h}d^{2h}(k^2+l^2)=\Box$, and conversely
$(k,l):=(d^hl',\,c^hk')$ satisfies~\eqref{eq:def}. Hence both formulations lead to the same curve $\Ccdh$.}
\end{rmk}

\begin{rmk}\label{rmk:ib}
{\rmksty Formally, $\Ccdh$ is obtained from the curve
$y^2=x^3+(a^{2h}+b^{2h})x^2+a^{2h}b^{2h}x$ studied in~\cite{HHZ} by replacing $(a^h,b^h)$ by $(c^h,{\rm i}\,d^h)$.
Notice also that $\Ccdh\cong\C_{(\lambda c)^h,(\lambda d)^h}$ over~$\Q$ for every positive integer~$\lambda$
(via $x\mapsto\lambda^{2h}x$, $y\mapsto\lambda^{3h}y$), and that $(c,d)$ is a {\hlppp} if and only if
$(\lambda c,\lambda d)$ is. Since the substitution $(m,n)\mapsto(m+n,m-n)$ transforms $(a,b,c)$ into $(2b,2a,2c)$,
the family of {\hls} of Case~1 coincides, up to a factor~$2$, with the family of Case~2.
For a fixed primitive triple, however, the {\hls} $(c,a)$ (with $c,a$ odd) and $(c,b)$ (with $b$ even)
are genuinely different (see the example above).}
\end{rmk}

\section{The torsion group of \texorpdfstring{$\boldsymbol{\Ccdh}$}{Gamma}}\label{sec:torsion}

We shall frequently use the following duplication identities.
Let $E\colon y^2=x(x+M)(x+N)$ with $M,N\in\Q\setminus\{0\}$, $M\neq N$.
If $P=(x_1,y_1)\in E(\Q)$ with $y_1\neq 0$, then we have
\begin{equation}\label{eq:dup}
	\begin{aligned}
	x_2&=\Bigl(\frac{x_1^2-MN}{2y_1}\Bigr)^2, \\
	x_2+M &=\Bigl(\frac{x_1^2+2Mx_1+MN}{2y_1}\Bigr)^2, \\
	x_2+N &=\Bigl(\frac{x_1^2+2Nx_1+MN}{2y_1}\Bigr)^2,
	\end{aligned}
\end{equation}
where $x_2$ denotes the $x$-coordinate of $[2]P:=P+P$. The first formula is the usual duplication formula (see, {\eg}, 
\cite[p.\,27]{SilvermanTate}),
and the other two follow from it by a direct computation using $y_1^2=x_1(x_1+M)(x_1+N)$.
For $\Ccdh$ we have $M=c^{2h}$ and $N=-d^{2h}$.

We shall also need the following arithmetic fact.

\begin{lem}\label{lem:nosquare}
Let $(c,d)$ be an {\hl} and $h$ a positive integer. Then $c^{2h}+d^{2h}$ is not a square.
\end{lem}

\begin{proof}
Assume that $c^{2h}+d^{2h}=z^2$ for some positive integer $z$, and let $g:=\gcd(c,d)$. Then
$(z/g^h)^2=(c/g)^{2h}+(d/g)^{2h}$ is an integer, hence so is $z/g^h$, and $(c/g,d/g)$ is again an {\hl}.
Thus we may assume that $\gcd(c,d)=1$. Since $c^2-d^2$ is a non-zero square, we have $c>d$.
For $h=1$, we would obtain that $c^4-d^4=(c^2+d^2)(c^2-d^2)$ is a non-zero square, contradicting
Fermat's Theorem that $x^4-y^4=z^2$ has no solution in positive integers.
For $h\geq 2$, the equation $c^{2h}+d^{2h}=z^2$ with $\gcd(c,d)=1$ is impossible
by~\cite[Main\;Theorem\,2]{Darmon_Merel}, applied with exponent $2h\geq 4$
(for even $h$ this follows already from Fermat's Theorem on $x^4+y^4=z^2$).
\end{proof}

\begin{prop}\label{prop:torsion}
Let $(c,d)$ be an {\hl} and $h$ a positive integer. Then the elliptic curve
$$\begin{aligned}
	\Ccdh: \ y^2 &= x(x+c^{2h})(x-d^{2h})
\end{aligned}$$
has torsion group $\tor 2 2$, consisting of the point $\mathscr{O}$ at infinity and the three points
$(0,0)$, $(-c^{2h},0)$, and $(d^{2h},0)$ of order~$2$.
\end{prop}

\begin{proof}
Write $C:=c^h$ and $A:=d^h$.
Since all three points of order~$2$ are rational, 
Mazur's Theorem~\cite{Mazur} implies that the
torsion group of $\Ccdh$ is $\tor 2{2N}$ with $N\in\{1,2,3,4\}$. Thus it is enough to show that
$\Ccdh$ has neither a point of order~$4$ nor a point of order~$3$.

{\it No point of order $4$:} If $Q$ is a point of order~$4$, then $T:=[2]Q$ is one of the three points of order~$2$,
and by~\eqref{eq:dup}, $x(T)$, $x(T)+C^2$, and $x(T)-A^2$ are squares of rational numbers.
For $T=(0,0)$ this would mean that $-A^2$ is a square, for $T=(-C^2,0)$ that $-C^2$ is a square, both of which are absurd.
For $T=(A^2,0)$ we would get $c^{2h}+d^{2h}=A^2+C^2=\Box$, which contradicts Lemma~\ref{lem:nosquare}.

{\it No point of order $3$:} Suppose that $P=(x_0,y_0)$ has order~$3$. Then $y_0\neq 0$ and $[2]P=-P$,
hence $x([2]P)=x_0$. By~\eqref{eq:dup} (with $M=C^2$ and $N=-A^2$), the rational numbers
$$\begin{aligned}
	U &:=\frac{x_0^2+C^2A^2}{2y_0} &\text{and} && V &:=\frac{x_0^2+2C^2x_0-C^2A^2}{2y_0}
\end{aligned}$$
satisfy $U^2=x_0$ and $V^2=x_0+C^2$, and $x_0-A^2=W^2$ for some $W\in\Q$. Since
$y_0^2=x_0(x_0+C^2)(x_0-A^2)=(UVW)^2\neq 0$, we may choose the sign of $W$ such that $y_0=UVW$. Then
$$\begin{aligned}
	U+V &=\frac{x_0(x_0+C^2)}{y_0}=\frac{U^2V^2}{UVW}=\frac{UV}{W},
\end{aligned}$$
so $W=\frac{UV}{U+V}$, and therefore
$$\begin{aligned}
	C^2 &=V^2-U^2 &\text{and} && A^2 &=U^2-W^2=\frac{U^3(U+2V)}{(U+V)^2}.
\end{aligned}$$
Let $t:=\frac{V+C}{U}$. Since $\frac{V+C}{U}\cdot\frac{V-C}{U}=1$, we have $t\neq 0$ and
$\frac VU=\frac12\bigl(t+t^{-1}\bigr)$, and hence
$$\begin{aligned}
	t(t^2+t+1) &= t^2\Bigl(1+\frac{2V}{U}\Bigr)=t^2\frac{U+2V}{U}=\Bigl(\frac{tA(U+V)}{U^2}\Bigr)^2.
\end{aligned}$$
Write $t= p/q$ with relatively prime integers $p,q$ and $q>0$. Then $pq\,(p^2+pq+q^2)$ is a square.
Since $p^2+pq+q^2>0$ and the three factors are pairwise relatively prime, we get $p=r^2$, $q=s^2$, and
$r^4+r^2s^2+s^4=\Box$ with $rs\neq 0$. This contradicts a theorem of Euler's 
which states that $x^4+x^2y^2+y^4=z^2$
has no solution in positive integers (see, {\eg}, 
\cite[Chap.\,4, Eqn.\,($7'$)]{Mordell}).
\end{proof}

\begin{rmk}
{\rmksty The argument for the non-existence of points of order~$3$ uses only that the curve has the form
$y^2=x(x+C^2)(x-A^2)$ with $C,A\in\Q\setminus\{0\}$. The arithmetic of the {\hl} $(c,d)$ is only needed
to exclude points of order~$4$. }
\end{rmk}

\section{Pytha\-potent {\hls} of degree~\texorpdfstring{$\boldsymbol{h}$}{h}}\label{sec:main}
We establish a correspondence between non-torsion rational points on $\Ccdh$ and pyth\-agorean pairs $(k,l)$ for which $(c^hk, d^hl)$ is an hl-pair, characterizing pyth\-apotent hl-pairs of degree~$h$ via the rank of $\Ccdh$ over~$\Q$.

\begin{thm}\label{thm:main}
An {\hl} $(c,d)$ is a {\hlppp} if and only if the elliptic curve $\Ccdh$ has positive rank over $\Q$.
\end{thm}

The proof is based on the following two lemmata, which give an explicit correspondence between
rational points on $\Ccdh$ and {\pps} $(k,l)$ such that $(c^hk,d^hl)$ is an {\hl}.

\begin{lem}\label{lem:point2pair}
Let $(c,d)$ be an {\hl}, let $h\geq 1$, and let $P=(x_1,y_1)$ be a rational point on $\Ccdh$ with $y_1\neq 0$. Put
$$\begin{aligned}
	k &:=x_1^2+c^{2h}d^{2h} & \text{and} && l &:=2c^hy_1.
\end{aligned}$$
Then
\begin{equation} \label{eq:kl}
	\begin{aligned}
	k^2+l^2 &=
	\bigl(x_1^2+2c^{2h}x_1-c^{2h}d^{2h}\bigr)^2,\\
	(c^hk)^2-(d^hl)^2 &=c^{2h}\bigl(x_1^2-2d^{2h}x_1-c^{2h}d^{2h}\bigr)^2\neq 0. 
	\end{aligned}
\end{equation}
Hence, after multiplying with a common denominator and taking absolute values, $(k,l)$ is a {\pp} such that
$(c^hk,d^hl)$ is an {\hl}. Moreover, $x([2]P)=\bigl(\frac{c^hk}{l}\bigr)^2$.
\end{lem}

\begin{proof}
Write $C:=c^h$ and $A:=d^h$, so that $y_1^2=x_1(x_1+C^2)(x_1-A^2)$. Clearly $k=x_1^2+C^2A^2>0$,
and $l\neq 0$ since $y_1\neq 0$. By the first identity in~\eqref{eq:dup} (with $M=C^2$ and $N=-A^2$)
we have $x([2]P)=\frac{k^2}{4y_1^2}=\bigl(\frac{Ck}{l}\bigr)^2$. Inserting this into the second and the
third identity in~\eqref{eq:dup} and multiplying by $4y_1^2$ and $4C^2y_1^2$, respectively, we obtain
the two identities in~\eqref{eq:kl}.

It remains to show that $(Ck)^2-(Al)^2\neq 0$. Assume to the contrary that $(Ck)^2=(Al)^2$, {\ie},
$\frac{k^2}{l^2}=\frac{A^2}{C^2}$. By the first identity in~\eqref{eq:kl}, $k^2+l^2$ is the square of a
rational number, and hence so is
$$\begin{aligned}
	\frac{A^2+C^2}{C^2} &=\frac{k^2}{l^2}+1=\frac{k^2+l^2}{l^2}.
\end{aligned}$$
Consequently, $c^{2h}+d^{2h}=A^2+C^2$ is the square of a rational number and hence, being an integer,
the square of an integer. This contradicts Lemma~\ref{lem:nosquare}.
\end{proof}

\begin{lem}\label{lem:pair2point}
Let $(c,d)$ be an {\hl} and $(k,l)$ a {\pp} such that $(c^hk,d^hl)$ is an {\hl}, say
$k^2+l^2=j^2$ and $(c^hk)^2-(d^hl)^2=w^2$ with $j,w>0$. Then
$$\begin{aligned}
	R &:=\Bigl(\frac{c^{2h}k^2}{l^2},\frac{c^{2h}kjw}{l^3}\Bigr)
\end{aligned}$$
is a rational point on $\Ccdh$ of infinite order.
\end{lem}

\begin{proof}
With $X:=\frac{c^{2h}k^2}{l^2}$ we have $X+c^{2h}=\frac{c^{2h}j^2}{l^2}$ and $X-d^{2h}=\frac{w^2}{l^2}$. Hence
$X(X+c^{2h})(X-d^{2h})=\bigl(\frac{c^{2h}kjw}{l^3}\bigr)^2$, {\ie}, $R\in\Ccdh(\Q)$. Since $k,j,w\neq 0$, the
$y$-coordinate of $R$ is not zero, so $R$ is not a point of order~$2$, and by Proposition~\ref{prop:torsion},
$R$ has infinite order.
\end{proof}

\begin{proof}[Proof of Theorem\;\ref{thm:main}]
If $(c,d)$ is a {\hlppp}, then Lemma~\ref{lem:pair2point} provides a rational point of infinite order on $\Ccdh$,
so the rank of $\Ccdh$ is positive. Conversely, assume that $P\in\Ccdh(\Q)$ has infinite order. For every integer $i\geq 1$,
the point $[i]P$ is not of order~$2$, so Lemma~\ref{lem:point2pair} applied to $[i]P$ gives a {\pp} $(k_i,l_i)$ such that
$(c^hk_i,d^hl_i)$ is an {\hl}, and $\bigl(\frac{c^hk_i}{l_i}\bigr)^2=x([2i]P)$. If $(k_i,l_i)$ and $(k_j,l_j)$ were
multiples of each other, then $x([2i]P)=x([2j]P)$, hence $[2i]P=\pm[2j]P$, which implies $i=j$. Thus the {\pps}
$(k_i,l_i)$, $i\geq 1$, are pairwise not multiples of each other, and at most one of them is a multiple of $(d,e)$.
\end{proof}

The proof shows in fact the following:

\begin{cor}
If $(c,d)$ is a {\hlppp}, then there are
infinitely many {\pps} $(k,l)$, not multiples of each other, such that
$(c^hk,d^hl)$ is an {\hl}.
\end{cor}

\begin{alg}\label{alg1}{\algsty
		The following algorithm constructs {\pps} $(k,l)$ 
		from rational points on $\Ccdh$ of infinite order,
		such that $(c^hk,d^hl)$ is an {\hl}.
		
		\begin{itemize}
			\item Let $P$ be a rational point on $\Ccdh$ of infinite 
			order and let $(x_i,y_i)$ be the coordinates of $[i]P$, 
			where $i$ is a positive integer. 
			\item Let $k_i$ and $l_i$ be relatively prime positive integers such that
			$$\begin{aligned}
			 \frac{k_i}{l_i} &=\Bigl|\frac{x_i^2+c^{2h}d^{2h}}{2c^hy_i}\Bigr|.
			\end{aligned}$$
		\end{itemize}
		Then $(k_i,l_i)$ is a {\pp} and $(c^hk_i,d^hl_i)$ is an {\hl}, and the pairs $(k_i,l_i)$ are not multiples
		of each other.
	}
\end{alg}

\begin{rmk}
{\rmksty Alternatively, one can proceed as in~\cite{HHZ}: if $x([2]P)=(g/f)^2$, then $\frac{k}{l}=\frac{g}{c^hf}$.
In the formulation of Remark~\ref{rmk:B}, the point $P=(x_1,y_1)$ yields the {\hl}
$(l',k')=(x_1^2+c^{2h}d^{2h},\,2d^hy_1)$, for which $(c^hk')^2+(d^hl')^2$ is a square.}
\end{rmk}

\section{Examples}\label{sec:examples}

We consider the pytha\-gorean triple $(a,b,c)=(3,4,5)$, {\ie}, $m=2$ and $n=1$.
Recall that a point $P=(x_1,y_1)$ yields the {\pp} $k^2+l^2=\Box$ with $k:l=(x_1^2+c^{2h}d^{2h}):(2c^hy_1)$ 
such that $(c^hk)^2-(d^hl)^2$ is a square, and, in the 
formulation of Remark~\ref{rmk:B}, the {\hl} $(l',k')$ with $l'^2-k'^2=\Box$ and $l':k'=(x_1^2+c^{2h}d^{2h}):(2d^hy_1)$, for which
$(c^hk')^2+(d^hl')^2$ is a square. In the cases where the latter pair is considerably smaller, we give both pairs.

\noindent{\bf Case 1: the {\hl} $(c,a)=(5,3)$.}
\begin{enumerate}
\item The curve $\C_{5,3}$ has rank~$1$ with 
generator $P=(-15,60)$. The point $P$ gives the trivial
pair $(k,l)=(3,4)=(a,b)$, whereas $[2]P=\bigl(\frac{225}{16},-\frac{3375}{64}\bigr)$ gives
$(k,l)=(481,600)$ with $481^2+600^2=769^2$ and
$$\begin{aligned}
	(5\cdot 481)^2-(3\cdot 600)^2 &=1595^2.
\end{aligned}$$
In particular, $(5,3)$ is a {\hlpppd 1}.
\item The curve $\C_{5^2,3^2}$ has rank~$1$ with generator $P=(-135,3780)$, which leads to
$(k,l)=(51,140)$ with $51^2+140^2=149^2$ and
$$\begin{aligned}
	(5^2\cdot 51)^2-(3^2\cdot 140)^2 &=195^2.
\end{aligned}$$
In the formulation of Remark~\ref{rmk:B}, the same point gives the {\hl} $(l',k')=(85,84)$ with 
$85^2-84^2=13^2$ and
$$\begin{aligned}
	(5^2\cdot 84)^2+(3^2\cdot 85)^2 &=2235^2.
\end{aligned}$$
In particular, $(5,3)$ is a {\hlpppd 2}  (compare Corollary~\ref{cor:h2}).
\item The curve $\C_{5^3,3^3}$ has rank~$1$ with generator $P=(-3600,432900)$, which leads to
$(k,l)=(9,40)$ with $9^2+40^2=41^2$ and
$$\begin{aligned}
	(5^3\cdot 9)^2-(3^3\cdot 40)^2 &=315^2.
\end{aligned}$$
In particular, $(5,3)$ is a {\hlpppd 3}.
\item The curve $\C_{5^4,3^4}$ has rank~$1$ with generator
$P=\bigl(\frac{1658842537500}{57410929},\frac{7157288054717943750}{435002609033}\bigr)$, which leads to
$k=594090016299$, $l=3596048594660$ with $k^2+l^2=3644791961501^2$ and
$$\begin{aligned}
	(5^4\cdot 594090016299)^2-(3^4\cdot 3596048594660)^2 &=230270140574595^2.
\end{aligned}$$
In particular, $(5,3)$ is a {\hlpppd 4}.
\item The curve $\C_{5^5,3^5}$ has rank~$1$ with generator
$P=\bigl(\frac{19775390625}{267289},\frac{14409551953125000}{138188413}\bigr)$, which leads to
$(k,l)=(1546371479601,\,1731239125000)$ with $k^2+l^2=2321304301649^2$ and
$$\begin{aligned}
	(5^5\cdot 1546371479601)^2-(3^5\cdot 1731239125000)^2 &=4814064171253125^2.
\end{aligned}$$
In the formulation of Remark~\ref{rmk:B}, the same point gives the smaller {\hl} 
$(l',k')=(19091005921, 1661989560)$ with $l'^2-k'^2=19018525121^2$ and
$$\begin{aligned}
	(5^5\cdot 1661989560)^2+(3^5\cdot 19091005921)^2 &=6963912904947^2.
\end{aligned}$$
In particular, $(5,3)$ is a {\hlpppd 5}.
\item For $h=6$ and $h=8$, the rank of $\C_{5^h,3^h}$ is~$0$. Hence, $(5,3)$ is not a {\hlppp} for these~$h$.
For $h=7$, Magma found the lower bound~$0$ and the upper bound~$1$ for the rank of  $\C_{5^7,3^7}$.
\end{enumerate}

\noindent{\bf Case 2: the {\hl} $(c,b)=(5,4)$.}
\begin{enumerate}\setcounter{enumi}{6}
\item The curve $\C_{5,4}$ has rank~$1$ with generator $P=(-20,60)$. The point $P$ gives the trivial
pair $(k,l)=(4,3)=(b,a)$, whereas $[2]P=\bigl(\frac{400}{9},-\frac{8000}{27}\bigr)$ gives
$(k,l)=(481,600)$ with $481^2+600^2=769^2$ and
$$\begin{aligned}
	(5\cdot 481)^2-(4\cdot 600)^2 &=155^2.
\end{aligned}$$
In particular, $(5,4)$ is a {\hlpppd 1}.
\item The curve $\C_{5^2,4^2}$ has rank~$0$. Hence, $(5,4)$ is not a {\hlpppd 2}.
As a matter of fact we would like to mention that 
for the same pythagorean triple $(3,4,5)$, the {\hl} $(5,3)$ is pytha\-potent of degree~$2$ while $(5,4)$ is not. 

\item The curve $\C_{5^3,4^3}$ has rank~$2$ with generators $P_1=(-12960,767520)$ and $P_2=(-6400,787200)$.
The point $P_1$ leads to $(k,l)=(272,225)$ with $272^2+225^2=353^2$ and
$$(5^3\cdot 272)^2-(4^3\cdot 225)^2=30800^2;$$
in the formulation of Remark~\ref{rmk:B}, $P_1$ gives the {\hl} $(l',k')=(85,36)$ with $85^2-36^2=77^2$ and
$$\begin{aligned}
	(5^3\cdot 36)^2+(4^3\cdot 85)^2 &=7060^2.
\end{aligned}$$
The point $P_2$
leads to $(k,l)=(8,15)$ with $8^2+15^2=17^2$ and
$$\begin{aligned}
	(5^3\cdot 8)^2-(4^3\cdot 15)^2 &=280^2.
\end{aligned}$$
In particular, $(5,4)$ is a {\hlpppd 3}.
\item The curve $\C_{5^4,4^4}$ has rank~$1$ with generator $P=(72200,14922600)$, which leads to
$(k,l)=(181252,109725)$ with $k^2+l^2=211877^2$ and
$$\begin{aligned}
	(5^4\cdot 181252)^2-(4^4\cdot 109725)^2 &=109744700^2.
\end{aligned}$$
In particular, $(5,4)$ is a {\hlpppd 4}.
\item The curve $\C_{5^5,4^5}$ has rank~$1$ with generator
$$P=\Bigl(\frac{106965627136000000000000}{35252663636483329},\frac{58126172986847820193382400000000000}{6618931792950754939769983}\Bigr),$$
which leads to
$$\begin{aligned}
k&=57563377462128924883445904754201,\\
l&=162466529187919179570759473250000
\end{aligned}$$
with $k^2+l^2=172362744034830586631245989777049^2$ and
\begin{multline*}
	(5^5\cdot 57563377462128924883445904754201)^2\\
	-(4^5\cdot 162466529187919179570759473250000)^2 \\
	={68419719323216368170528698116878125}^2.
\end{multline*}
In particular, $(5,4)$ is a {\hlpppd 5}.
\item The curve $\C_{5^6,4^6}$ has rank~$0$.
The curve $\C_{5^7,4^7}$ has positive rank: the point
$$\begin{aligned}
	P &=\Bigl(\frac{163840000000000}{210681},\frac{5048068833280000000000}{96702579}\Bigr)
\end{aligned}$$
has infinite order and leads to
$(k,l)=(9199112,\,33449625)$ with $k^2+l^2=34691513^2$ and
$$\begin{aligned}
	(5^7\cdot 9199112)^2-(4^7\cdot 33449625)^2 &=464925233000^2.
\end{aligned}$$
In particular, $(5,4)$ is a {\hlpppd 7}.
\end{enumerate}

\section{The degrees~\texorpdfstring{$\boldsymbol{h=1,2,3}$}{h=1,2,3}}\label{sec:small}

We now show that every {\hl} is a {\hlpppd 1} as well
as a {\hlpppd 3}, and we give sufficient conditions for degree~$2$.
In each case we first formulate the result for an arbitrary {\hl} $(c,d)$ with $e=\sqrt{c^2-d^2}$, and then make it explicit
for the two cases $(c,d)=(c,a)$ and $(c,d)=(c,b)$, where $a=m^2-n^2$, $b=2mn$, $c=m^2+n^2$.

\begin{cor}\label{cor:h1}
Every {\hl} $(c,d)$ is a {\hlpppd 1}. More precisely, the point $P_0=(c^2 d^2/e^2, c^3 d^3/e^3)$ on $\Ccd$ is a point of infinite order
and yields the {\pp} $(k,l)$ with
$$\begin{aligned}
	k &=d^4+d^2e^2+e^4, & l &=2c^2de,
\end{aligned}$$
for which
$$\begin{aligned}
	k^2+l^2 &=\bigl(c^4+d^2e^2\bigr)^2 & \text{and} &&
	(ck)^2-(dl)^2 &=\bigl(c(e^4+d^2e^2-d^4)\bigr)^2.
\end{aligned}$$
\end{cor}

\begin{proof}
Using $c^2=d^2+e^2$, one checks that $P_0$ lies on $\Ccd$, and since its $y$-coordinate is non-zero,
$P_0$ has infinite order by Proposition~\ref{prop:torsion}. In fact, $P_0=[2]Q$ for the point
$Q=(-cd,-cde)\in\Ccd(\Q)$, which by Lemma~\ref{lem:point2pair} yields the trivial
{\pp} $(k,l)\sim(d,e)$ (compare the examples in Section~\ref{sec:examples}).
By Lemma~\ref{lem:point2pair}, $P_0$ gives 
$$\begin{aligned}
	k' &=c^2d^2+\frac{c^4d^4}{e^4}
	&\text{and}&&
	l' &=\frac{2c^4d^3}{e^3},
\end{aligned}$$
and multiplying by $\frac{e^4}{c^2d^2}$, we obtain the
{\pp} $(k,l)$ with 
$$\begin{aligned}
	k &=d^4+d^2e^2+e^4, & l &=2c^2de,
\end{aligned}$$
where we used the identity $c^2d^2=d^4+d^2e^2$.

The two identities for $k^2+l^2$ and $(ck)^2-(dl)^2$
follow from~\eqref{eq:kl} and $c^2=d^2+e^2$. Finally, $(k,l)$ is not a 
multiple of the trivial {\pp} $(d,e)$,
since otherwise, $e^4=c^2d^2$, {\ie}, $(c/d)^2-(c/d)-1=0$, 
which is impossible for rational $c/d$.
\end{proof}

In terms of $m$ and $n$ we obtain in both cases ({\ie}, $d=m^2-n^2$
and $d=2mn$), the same {\pp}
$$\begin{aligned}
	k &=m^8+14m^4n^4+n^8, & l &=4mn(m^4-n^4)(m^2+n^2),
\end{aligned}$$
with 
$$\begin{aligned}
	k^2+l^2 &=\bigl(m^8+8m^6n^2-2m^4n^4+8m^2n^6+n^8\bigr)^2.
\end{aligned}$$
Finally, for $(ck)^2-(dl)^2=\bigl((m^2+n^2)k\bigr)^2-\bigl((m^2-n^2)l\bigr)^2$ 
we obtain the square
$$\begin{aligned}
	&\bigl((m^2+n^2)(m^4-2m^3n-2m^2n^2-2mn^3+n^4)(m^4+2m^3n-2m^2n^2+2mn^3+n^4)\bigr)^2,
\end{aligned}$$
and for $(ck)^2-(dl)^2=\bigl((m^2+n^2)k\bigr)^2-\bigl(2mnl\bigr)^2$
we obtain the square
$$\begin{aligned}
	&\bigl((m^2+n^2)(m^4-4m^2n^2-n^4)(m^4+4m^2n^2-n^4)\bigr)^2.
\end{aligned}$$

For $h=3$, we get a similar result as for $h=1$:

\begin{cor}\label{cor:h3}
Every {\hl} $(c,d)$ is a {\hlpppd 3}. More precisely, the point 
$P_0=(c^4d^2,c^4d^2e(c^2 + d^2))$ on $\C_{c^3,d^3}$ has infinite order and
it yields the {\pp} $(k,l)$ with
$$\begin{aligned}
	k &=d^2 & \text{and} && l &=2 c e,
\end{aligned}$$
for which we have
$$
\begin{aligned}
	k^2+l^2 &=\bigl(2 c^2-d^2\bigr)^2, \\
	(c^3k)^2-(d^3 l)^2 &=\bigl(d^2 (c^3 - 2 c d^2)\bigr)^2.
\end{aligned}$$
\end{cor}

\begin{proof}
Using $e^2=c^2-d^2$, one checks that $P_0$ lies on $\C_{c^3,d^3}$, and since its $y$-coordinate is non-zero,
$P_0$ has infinite order by Proposition~\ref{prop:torsion}.
By Lemma~\ref{lem:point2pair}, $P_0$ gives 
$$\begin{aligned}
	k' &=c^6 d^4 (c^2 + d^2) & \text{and} &&	l' &=2 c^7 d^2 (c^2 + d^2) e,
\end{aligned}$$
and dividing by $c^6d^2(c^2+d^2)$, we obtain the
{\pp} $(k,l)$ with 
$$\begin{aligned}
	k &=d^2, & l &=2ce,
\end{aligned}$$
which is obviously not a multiple of the {\pp} $(d,e)$.
The other two identities follow easily from $e^2=c^2-d^2$.
\end{proof}


In terms of $m$ and $n$, for $d=m^2-n^2$ we obtain the {\pp}
$$\begin{aligned}
	k &=(m^2 - n^2)^2, & l &=4 m n (m^2 + n^2),
\end{aligned}$$
with 
$$\begin{aligned}
	k^2+l^2 &=(m^4 + 6 m^2 n^2 + n^4)^2,
\end{aligned}$$
and for $d=2mn$, after dividing by $2$ we obtain the {\pp}
$$\begin{aligned}
	k &=2 m^2 n^2, & l &=(m^2 - n^2) (m^2 + n^2),
\end{aligned}$$
with 
$$\begin{aligned}
	k^2+l^2 &= (m^4 + n^4)^2.
\end{aligned}$$
Finally, for $(c^3k)^2-(d^3l)^2=\bigl((m^2+n^2)^3k\bigr)^2-\bigl((m^2-n^2)^3l\bigr)^2$ 
we obtain the square
$$\begin{aligned}
	& \bigl((m^2+n^2)(m^2-n^2)^2(m^4-6m^2n^2+n^4)\bigr)^2,
\end{aligned}$$
and for $(c^3k)^2-(d^3l)^2=\bigl((m^2+n^2)^3k\bigr)^2-\bigl((2mn)^3l\bigr)^2$
we obtain the square
$$\begin{aligned}
	& \bigl(2m^2n^2(m^2+n^2)(m^4-6m^2n^2+n^4)\bigr)^2.
\end{aligned}$$

Now, we consider also the case when $h=2$:

\begin{cor}\label{cor:h2}
Let $(c,d)$ be an {\hl} with $e=\sqrt{c^2-d^2}$, and assume that for some $\varepsilon\in\{1,-1\}$
$$\begin{aligned}
	c^2+\varepsilon cd+d^2 &=q^2 & \text{for some positive integer $q$.}
\end{aligned}$$
Then $(c,d)$ is a {\hlpppd 2}. More precisely, $P=(-\varepsilon cd^3, cd^3eq)$ is a point of infinite order
on $\C_{c^2,d^2}$ and yields the {\pp}
$$\begin{aligned}
	k &=d(c^2+d^2), & l &=2ceq,
\end{aligned}$$
for which
$$\begin{aligned}
	k^2+l^2 &=\bigl(2c^3+\varepsilon(c^2d-d^3)\bigr)^2 &\text{and} && (c^2k)^2-(d^2l)^2 &=\bigl(cd\,(c^3-cd^2-2\varepsilon d^3)\bigr)^2.
\end{aligned}$$
\end{cor}

\begin{proof}
We have $(-\varepsilon cd^3)(c^4-\varepsilon cd^3)(-\varepsilon cd^3-d^4)=c^2d^6(c^3-\varepsilon d^3)(c+\varepsilon d)
=c^2d^6(c^2-d^2)(c^2+\varepsilon cd+d^2)=(cd^3eq)^2$, so $P$ is a point on $\C_{c^2,d^2}$ of infinite order
by Proposition~\ref{prop:torsion}, and the claim follows from Theorem~\ref{thm:main}. By Lemma~\ref{lem:point2pair},
$P$ yields $k=c^2d^6+c^4d^4=c^2d^4(c^2+d^2)$ and $l=2c^3d^3eq$, {\ie}, $(k,l)\sim(d(c^2+d^2),\,2ceq)$, and the
two identities follow from $e^2=c^2-d^2$ and $q^2=c^2+\varepsilon cd+d^2$.
\end{proof}

In terms of $m$ and $n$, the conditions of Corollary~\ref{cor:h2} read as follows:
\begin{itemize}
\item[] {\bf Case 1:} The {\hl} $(c,a)$ is a {\hlpppd 2} if 
at least one of the following two conditions is satisfied:
$$\begin{aligned}
	\text{(i)} \ 3m^4+n^4 &=\Box \\ \text{(ii)} \ m^4+3n^4 &=\Box
\end{aligned}$$

If $q^2$ is the respective square, then $(k,l)=\bigl((m^2-n^2)(m^4+n^4),\;2mn(m^2+n^2)\,q\bigr)$.
\item[] {\bf Case 2:} The {\hl} $(c,b)$ is a {\hlpppd 2}  
if at least one of the following two conditions is satisfied:
	$$\begin{aligned}
		\text{(i)} \ m^4+2m^3n+6m^2n^2+2mn^3+n^4 &=\Box\\
		\text{(ii)} \ m^4-2m^3n+6m^2n^2-2mn^3+n^4 &=\Box
	\end{aligned}$$
If $q^2$ is the respective square, then 
$(k,l)=\bigl(mn(m^4+6m^2n^2+n^4), (m^4-n^4)\,q\bigr)$.
\end{itemize}
Each of these four quartic conditions defines a curve of genus~$1$ with a rational point which is
birational to the elliptic curve $y^2=x^3-12x$ of rank~$1$; hence each condition is satisfied by infinitely many
pairs~$(m,n)$. For example, $(m,n)=(2,1)$ satisfies Case~1\,(i) and gives the {\hl} $(5,3)$
(see Section~\ref{sec:examples}), and $(m,n)=(11,3)$ satisfies Case~1\,(ii) and gives the {\hl} $(130,112)$, {\ie}, $(65,56)$.
In Case~2, $(m,n)=(3,1)$ satisfies~(i) and gives $(10,6)$, {\ie}, again $(5,3)$, and $(m,n)=(7,4)$ satisfies~(ii) and
gives $(65,56)$ ({\sl cf}. Remark~\ref{rmk:ib}).

On the other hand, the {\hl} $(5,4)$ satisfies neither condition, in accordance with the fact that
$\C_{5^2,4^2}$ has rank~$0$.

We conclude the paper with two open problems.
\begin{ques}
{\rm Take an arbitrary {\hl} $(c,d)$. Is there an $h \geq 4$ such that $(c, d)$ is a pytha\-potent {\hl} of degree~$h$?}
\end{ques}
\begin{ques}
{\rm Take an arbitrary $h \geq 4$. Is there a {\hl} 
$(c, d)$ which is a pytha\-potent {\hl} of degree~$h$? 
Or equivalently: Is there a {\hl} $(c,d)$ 
such that $\Gamma_{c^h,d^h}$ has
positive rank over $\Q$?}
\end{ques}

The restriction to $h\geq 4$ in both questions is natural: by Corollaries~\ref{cor:h1} and~\ref{cor:h3},
every {\hl} is a {\hlpppd 1} and a {\hlpppd 3}, whereas for $h=2$ the answer depends on the {\hl};
for example, $(5,3)$ is a {\hlpppd 2}, but $(5,4)$ is not (see Corollary~\ref{cor:h2} and Section~\ref{sec:examples}).
The examples in Section~\ref{sec:examples} answer the second question for $h=4,5,7$;
for $h=6$ and $h\geq 8$ it remains open.

\subsection*{Statements and Declarations}
The authors have no competing interests to declare that are relevant to the content of this article.
No funding was received for conducting this study.

\subsection*{Declaration of generative AI and AI-assisted technologies in the manuscript preparation process}
During the preparation of this work, the authors used Claude Opus 5.5 to improve the language and to check correctness of proofs. 
The authors take full responsibility for the content of the published article. 

\bibliographystyle{plain}

\end{document}